\documentclass[11pt,a4paper,reqno]{amsart}
\usepackage{a4wide,amsfonts,amssymb,mathptmx}
\usepackage{hyperref}

\providecommand{\BBb}[1]{{\mathbb{#1}}}
\providecommand{\cal}[1]{{\mathcal{#1}}}   

\newcommand{\B}{{\BBb B}}

\newcommand{\lap}{\operatorname{\Delta}}
\newcommand{\mlap}{-\!\operatorname{\Delta}}
\newcommand{\N}{\BBb N}
\newcommand{\R}{{\BBb{R}}}
\newcommand{\Rn}{{\BBb{R}^{n}}}
\newcommand{\Span}{\operatorname{span}}

\newtheorem{thm}{Theorem}
 
\theoremstyle{remark}
\newtheorem{rem}[thm]{Remark}

\newcommand{\vvvert}{{|\hspace{-1.6pt}|\hspace{-1.6pt}|}}

\newcommand{\ip}[2]{\ensuremath{({#1}\,|\,{#2})}} 

\def\Xint#1{\mathchoice
{\XXint\displaystyle\textstyle{#1}}%
{\XXint\textstyle\scriptstyle{#1}}%
{\XXint\scriptstyle\scriptscriptstyle{#1}}%
{\XXint\scriptscriptstyle\scriptscriptstyle{#1}}%
\!\int}
\def\XXint#1#2#3{{\setbox0=\hbox{$#1{#2#3}{\int}$}
\vcenter{\hbox{$#2#3$}}\kern-.5\wd0}}

\def\dashint{\Xint-}

\begin{document}
\title[On final value problems and well-posedness]{%
On parabolic final value problems and well-posedness}
\author[Christensen and Johnsen]{Ann-Eva Christensen{{\ensuremath{^{\upn{a}\,1}}}} %
and Jon Johnsen{{\ensuremath{^{\upn{a}}}}}}
\address{$^1$Present address: Unit of Epidemiology and Biostatistics, Aalborg University Hospital, Hobrovej 18--22,  DK-9000 Aalborg, Denmark}
\email{anneva.ch@gmail.com}
\address{$^\upn{A}$Department of Mathematical Sciences, Aalborg University, Skjernvej 4A, DK-9220 Aalborg {\O}st, Denmark}
\email{{\rm(corresponding author)} jjohnsen@math.aau.dk}
\thanks{The second author is supported by the Danish 
Research Council, Natural Sciences grant no.~4181-00042.
\\[9\jot]{\tt Appeared online in C. R. Acad.\ Sci. Paris, Ser. I. (2018); DOI: 10.1016/j.crma.2018.01.019\ }}

\begin{abstract} We prove that a large class of parabolic final  value problems is well posed.  
  This results via explicit Hilbert spaces that characterise the data
  yielding existence, uniqueness and stability of solutions.  
  This data space is the graph normed domain of an unbounded operator, 
  which represents a new compatibility condition pertinent for final value problems.  
  The framework is evolution equations for
  Lax--Milgram operators in vector distribution spaces.
  The final value heat equation on a smooth open set is
  also covered, and for non-zero Dirichlet data a non-trivial extension of the
  compatibility condition is obtained by addition of an improper Bochner integral.
\end{abstract}
\maketitle
\section{Introduction}\label{intro-sect}
\thispagestyle{empty}

Well-posedness of final value problems for a large class of parabolic differential equations is described here.
That is, for suitable spaces $X$, $Y$ specified below, they have \emph{existence, uniqueness} and \emph{stability} of solutions $u\in X$ for given data $(f,g,u_T)\in Y$.
This should provide a basic clarification of a type of problems, which hitherto has been insufficiently understood.

As a first example, we 
characterise the functions $u(t,x)$ that, in a $C^\infty$-smooth bounded open set
$\Omega\subset\Rn$ with boundary $\partial\Omega$, 
satisfy the following equations that constitute the final value problem for the heat equation ($\lap=\partial_{x_1}^2+\dots+\partial_{x_n}^2$ denotes the Laplacian):
\begin{equation}  \label{heat-intro}
\left.
\begin{aligned}
  \partial_tu(t,x)-\lap u(t,x)&=f(t,x) &&\quad\text{for $t\in\,]0,T[\,$,  $x\in\Omega$,}
\\
   u(t,x)&=g(t,x) &&\quad\text{for $t\in\,]0,T[\,$, $x\in\partial\Omega$,}
\\
  u(T,x)&=u_T(x) &&\quad\text{for $x\in\Omega$}.
\end{aligned}
\right\}
\end{equation}
Hereby $(f,g,u_T)$ are the given data of the problem.

In case $f=0$, $g=0$ the first two lines of \eqref{heat-intro} are satisfied by 
$u(t,x)=e^{(T-t)\lambda}v(x)$  for all $t\in\R$, if $v(x)$ is an eigenfunction of the Dirichlet realization
$\mlap_D$ with eigenvalue $\lambda$. 

Thus the homogeneous final value problem
\eqref{heat-intro} has the above $u$ as a \emph{basic} solution if, coincidentally,
the final data $u_T$ equals the eigenfunction $v$.
Our construction includes the set $\cal B$ of such basic solutions $u$, its linear hull
$\cal E=\Span\cal B$ and a certain completion $\overline{\cal E}$. 

Using the eigenvalues $0<\lambda_1\le\lambda_2\le\dots$ and the associated $L_2(\Omega)$-orthonormal
basis $e_1, e_2,\dots$ of eigenfunctions  of $\mlap_D$, the space $\cal E$ (that corresponds to data $u_T\in\Span( e_j)$) clearly consists of solutions $u$ being \emph{finite} sums
\begin{equation}
  \label{basic-id}
  u(t,x)=\textstyle{\sum_j\;} e^{(T-t)\lambda_j}\ip{u_T}{e_j}e_j(x).
\end{equation}
So at $t=0$ there is, by the finiteness, a vector $u(0,x)$ in $L_2(\Omega)$ fulfilling
\begin{equation}
  \|u(0,\cdot)\|^2 =\textstyle{\sum_j\;} e^{2T\lambda_j}|\ip{u_T}{e_j}| ^2 <\infty.
  \label{basic-ineq}
\end{equation}

When summation is extended to all $j\in\N$, condition \eqref{basic-ineq} becomes very strong, as it is only satisfied for special $u_T$: 
by Weyl's law $\lambda_j= \cal O(j^{2/n})$, 
so a single term in \eqref{basic-ineq} yields $|\ip{u_T}{e_j}|\le c \exp({-Tj^{2/n}})$; whence the $L_2$-coordinates of such $u_T$ decay rapidly for $j\to\infty$.
This has been known since the 1950's; cf.\ the work of John~\cite{John55}
and Miranker~\cite{Miranker61}. 

More recently e.g.\ Isakov~\cite{Isa98} emphasized the observation, made already in
\cite{Miranker61}, that \eqref{basic-id} gives rise to an \emph{instability}: the sequence of data
data $u_{T,k}=e_k$ has length $1$ for all $k$, but \eqref{basic-id}  gives
$\|u_k(0,\cdot)\|=\|e^{T\lambda_k}e_k\|=e^{T\lambda_k}\nearrow\infty$ for $k\to\infty$. 
Thus \eqref{heat-intro} is not well-posed in $L_2(\Omega)$. 
 
 In general this instability shows that the $L_2$-norm is an insensitive choice.
To obtain well-adapted spaces for \eqref{heat-intro} with $f=0$, $g=0$, one could
depart from \eqref{basic-ineq}. Indeed, along with the solution space
$\cal E$, a norm on the final data $u_T\in\Span(e_j)$ can be \emph{defined}
by \eqref{basic-ineq}; and 
$\vvvert u_T \vvvert=(\sum_{j=1}^\infty  e^{2T\lambda_j}|\ip{u_T}{e_j}|^2)^{1/2}$ can be used as norm on the $u_T$
that correspond to solutions $u$ in the above completion $\overline{\cal E}$. This would give well-posedness of the 
homogeneous version of \eqref{heat-intro} with $u\in\overline{\cal E}$. (Cf.\ \cite{ChJo17}.)

But we have first of all replaced 
specific eigenvalue distributions by using sesqui-linear forms, cf.\ Lax--Milgram's lemma,
which allowed us to cover general elliptic operators $A$.

Secondly the \emph{fully} inhomogeneous problem \eqref{heat-intro} is covered. 
Here it does not suffice  to choose the norm on the data $(f,g,u_T)$
suitably (cf.~$\vvvert u_T\vvvert$), for 
one has to \emph{restrict} $(f,g,u_T)$ to a subspace first by imposing certain \emph{compatibility conditions}. These have long been  known for parabolic problems, but they
have a new form for final value problems.

\section{The abstract final value problem}
Our main analysis concerns a (possibly non-selfadjoint)
Lax--Milgram operator $A$ defined in $H$ from a bounded $V$-elliptic sesquilinear form $a(\cdot,\cdot)$ 
in a Gelfand triple, i.e.\
densely injected Hilbert spaces $V\hookrightarrow H\hookrightarrow V^*$ 
with norms $\|\cdot\|$, $|\cdot|$ and $\|\cdot\|_*$. 

In this set-up, we consider the following general final value problem: given data  
$ f\in L_2(0,T; V^*)$, $u_T\in H$, determine the vector distributions 
 $u\in\cal D'(0,T;V)$ fulfilling
\begin{equation}
  \label{eq:fvA-intro}
  \left.
  \begin{aligned}
  \partial_tu +Au &=f  &&\quad \text{in $\cal D'(0,T;V^*)$,}
\\
  u(T)&=u_T &&\quad\text{in $H$}.
\end{aligned}
\right\}
\end{equation}
A wealth of parabolic Cauchy problems with homogeneous boundary 
conditions have been efficiently treated using such triples $(H,V,a)$ and the $\cal D'(0,T;V^*)$
framework in \eqref{eq:fvA-intro}; cf.\ works of Lions and Magenes~\cite{LiMa72}, Tanabe~\cite{Tan79}, Temam~\cite{Tm}, Amann~\cite{Ama95}.
Also recently e.g.\ Almog, Grebenkov, Helffer, Henry  studied
variants of the complex Airy operator via such triples \cite{AlHe15,GrHelHen17,GrHe17}, and our results should at least extend to final value problems for those of their realisations that have non-empty spectrum.

For the corresponding Cauchy problem we recall that when solving $u'+Au=f$
so that  $u(0)=u_0$ in $H$, for $f\in L_2(0,T;V^*)$, there is a
unique solution $u$ in the Banach space 
\begin{equation}
  \begin{split}
  X:=&L_2(0,T;V)\bigcap C([0,T];H) \bigcap H^1(0,T;V^*)
\\
  \|u\|_X=&\big(\int_0^T (\|u(t)\|^2   +\|u'(t)\|_*^2)\,dt+\sup_{0\le t\le T}|u(t)|^2\big)^{1/2}.   
  \end{split}
  \label{eq:X-intro}
\end{equation}
For \eqref{eq:fvA-intro} it would therefore be natural to expect solutions $u$ in the same space $X$.
This is correct, but only when the data $(f,u_T)$ satisfy substantial further conditions.

To state these, we utilise that $-A$ generates an \emph{analytic} semigroup $e^{-tA}$ 
in $\B(H)$ and $\B(V^*)$, and that $e^{-tA}$ consequently is \emph{invertible} in the
class of closed operators on $H$, resp.\ $V^*$; cf.\ Proposition~2.2 in \cite{ChJo17}.
Consistently with the case when $A$ generates a group, we set
\begin{equation}
  \label{eq:inverse-intro}
  (e^{-tA})^{-1}=e^{tA}.
\end{equation}
Its domain $D(e^{tA})=R(e^{-tA})$ is the Hilbert space normed by $\|u\|=(|u|^2+|e^{tA}u|^2)^{1/2}$.
In the common case $A$ has non-empty spectrum, $\sigma(A)\ne\emptyset$, there is a chain of strict inclusions
\begin{equation} \label{dom-intro}
  D(e^{t'A})\subsetneq D(e^{t A})\subsetneq H \qquad\text{ for  $0<t<t'$}.
\end{equation}
At the final time $t=T$ these domains enter the well-posedness result below, where for breviety 
$y_f$ will denote the full yield of the source term $f$ on the system, namely
\begin{equation} \label{yf-intro}
  y_f= \int_0^T e^{-(T-s)A}f(s)\,ds.
\end{equation}
The map $f\mapsto y_f$ takes values in $H$, and it  is a continuous surjection $y_f\colon L_2(0,T;V^*)\to H$.

\begin{thm} \label{intro-thm}
  The final value problem \eqref{eq:fvA-intro} has a solution
  $u$ in the space $X$ in \eqref{eq:X-intro} if and only if the data
  $(f,u_T)$ belong to the subspace $Y$ of $L_2(0,T; V^*)\oplus H$ defined
  by the condition  
  \begin{equation}
    \label{eq:cc-intro}
    u_T-y_f \ \in\  D(e^{TA}).
  \end{equation}
In the affirmative case,  the solution $u$ is unique in $X$, and it depends continuously on the data $(f,u_T)$ in $Y$,
that is $\|u\|_X\le c\|(f,u_T)\|_Y$,
when $Y$ is given the graph norm
\begin{equation}
  \label{eq:Y-intro}
  \|(f,u_T)\|_Y=
  \Big(|u_T|^2+\int_0^T\|f(t)\|_*^2\,dt+\big|e^{TA}(u_T-y_f)\big|^2
  \Big)^{1/2}.
\end{equation}
\end{thm}

Condition \eqref{eq:cc-intro} is seemingly a fundamental novelty for the final value problem \eqref{eq:fvA-intro}.
As for \eqref{eq:Y-intro}, it is the graph norm of $(f,u_T)\mapsto e^{TA}(u_T-y_f)$,
which for $\Phi(f,u_T)=u_T-y_f$ is the unbounded operator $e^{TA}\circ\Phi$ from $L_2(0,T; V^*)\oplus H$ to $H$.  

In fact, $e^{TA}\Phi$ is central to a rigorous treatment of \eqref{eq:fvA-intro}, for \eqref{eq:cc-intro} means that $e^{TA}\Phi$ must be defined at
$(f,u_T)$; i.e.\ the data space $Y$ is its domain.
So since $e^{TA}\Phi$ is a closed operator, $Y$ is a Hilbert space, which 
 by \eqref{eq:Y-intro} is embedded into
$L_2(0,T; V^*)\oplus H$.

As an inconvenient aspect, the presence of $e^{-(T-t)A}$ and the integration over $[0,T]$ make
\eqref{eq:cc-intro} \emph{non-local} in space and time \,---\,exacerbated by use of the abstract domain $D(e^{TA})$, which for larger $T$ 
gives increasingly stricter conditions; cf.\ \eqref{dom-intro}. 

We regard \eqref{eq:cc-intro} as a \emph{compatibility} condition on
the data $(f,u_T)$, and thus we generalise the notion.
Grubb and Solonnikov~\cite{GrSo90} made a systematic treatment
of \emph{initial}-boundary problems of parabolic equations
with compatibility conditions, which are necessary and sufficient for
well-posedness in full scales of anisotropic $L_2$-Sobolev spaces\,---\,whereby compatibility conditions are decisive for the solution's regularity.
In comparison \eqref{eq:cc-intro} is crucial for the \emph{existence} question; cf.\ Theorem~\ref{intro-thm}.

\begin{rem}
Previously uniqueness was observed by Amann~\cite[V.2.5.2]{Ama95} in a $t$-dependent set-up. 
However, the injectivity of $u(0)\mapsto u(T)$ was shown much earlier in a set-up with $t$-dependent sesquilinear forms by Lions and Malgrange~\cite{LiMl60}.
\end{rem}
\begin{rem}
Showalter~\cite{Sho74} attempted to characterise the possible
$u_T$ in terms of Yosida approximations for $f=0$ and $A$ having half-angle $\pi/4$. As an ingredient the 
invertibility of analytic semigroups was claimed by Showalter for such $A$, but his proof was flawed as
$A$ can have semi-angle $\pi/4$ even if $A^2$ is not accretive; cf.\ our example in Remark~3.15 of \cite{ChJo17}.
\end{rem}

Theorem~\ref{intro-thm} is proved by considering the full set of solutions to the differential equation $u'+Au=f$. As indicated in \eqref{eq:X-intro}, for fixed $f\in L_2(0,T;V^*)$ the solutions in $X$ are parametrised by the initial state $u(0)\in H$; and they are also in this set-up necessarily given by the variation of constants formula for the analytic semigroup $e^{-tA}$ in $V^*$, 
\begin{equation}
  \label{eq:bijection-intro}
  u(t)=e^{-tA}u(0)+ \int_0^t e^{-(t-s)A}f(s)\,ds.
\end{equation}
For $t=T$ this yields a \emph{bijective correspondence} $u(0)\longleftrightarrow u(T)$ between the initial
and terminal states---for due to the invertibility of  $e^{-TA}$, cf.\ \eqref{eq:inverse-intro}, one can isolate $u(0)$ here.
Moreover, \eqref{eq:bijection-intro} also yields necessity of \eqref{eq:cc-intro} at once, as the
difference $u_T-y_f$ in \eqref{eq:cc-intro}  must be equal to $e^{-TA}u(0)$, which clearly belongs to the domain $D(e^{TA})$.  

Moreover, $u(T)$ consists of two 
radically different parts, cf.\  \eqref{eq:bijection-intro}, even when $A$ is `nice': 

First, $e^{-tA}u(0)$ solves the equation for $f=0$, and for $u(0)\ne0$ we obtained in \cite{ChJo17}
the precise property in non-selfadjoint dynamics that the ``height'' function $h(t)$ is \emph{strictly convex}. Hereby
\begin{equation}
  h(t)= |e^{-tA}u(0)|.
\end{equation}
This results from the injectivity of $e^{-tA}$ when
 $A$ is normal, or belongs to the class of hyponormal operators studied by Janas~\cite{Jan94}, 
or in case $A^2$ is accretive\,---\,so for such $A$ the complex eigenvalues (if any) cannot give oscillations in the size of $e^{-tA}u(0)$, cf.\ the strict convexity. This stiffness from the strict convexity is consistent with the fact for analytic semigroups  that $u(T)=e^{-TA}u(0)$ is
confined to the dense, but very small space $\bigcap_{n\in\N} D(A^n)$.

In addition $h(t)$ is strictly decreasing with $h'(0)\le-m(A)$, where $m(A)$ denotes the lower bound; i.e.\ the 
short-time behaviour is governed by the numerical range $\nu(A)$ also for such $A$.

Secondly, for $u(0)=0$ the equation is solved by the integral in   \eqref{eq:bijection-intro}, 
which has rather different properties. Its final value $y_f\colon L_2(0,T;V^*)\to H$ is surjective, so $y_f$ 
can be \emph{anywhere} in $H$.
This was shown with a kind of control-theoretic argument in \cite{ChJo17} for the case that $A=A^*$ with $A^{-1}$ compact;
and for general $A$ by using the Closed Range Theorem.

Thus the possible final data $u_T$ are a sum of an arbitrary $y_f\in H$ and a
term $e^{-TA}u(0)$ of great stiffness, so that $u_T$ can be prescribed anywhere in the affine space
$y_f+D(e^{TA})$. As $D(e^{TA})$  is dense in $H$, and in general there hardly is any control over the direction of $y_f$ (if non-zero), it is not feasible to specify $u_T$ a priori in other spaces than $H$. Instead it is by the condition $u_T-y_f\in D(e^{TA})$ that the $u_T$ and $f$ are properly controlled.

\section{The inhomogeneous heat problem}
For general data $(f,g,u_T)$ in \eqref{heat-intro}, the results in Theorem~1 are applied with
$A=\mlap_{D}$. The results are
analogous, but less simple to prove and state.

First of all, even though it is a linear problem,  
the compatibility condition \eqref{eq:cc-intro}
\emph{destroys} the old trick of reducing to boundary data $g=0$, 
for when $w\in H^1$ fulfils $w=g\ne0$ on the curved boundary $\,]0,T[\,\times\partial\Omega$, then $w$
\emph{lacks} the regularity needed to test condition \eqref{eq:cc-intro} on the resulting data 
$(\tilde f,0,\tilde u_T)$ of the reduced problem.

Secondly, it therefore takes an effort to show that when the boundary data $g\ne0$, then they \emph{do} give rise to a correction term $z_g$. This means that
condition \eqref{eq:cc-intro} is replaced by 
\begin{equation}
  u_T-y_f+z_g\in D(e^{-T\lap_{D}}).  
\end{equation}

Thirdly, because of the low reqularity, it requires some  technical diligence to show that, despite the singularity 
present in $\lap e^{(T-s) \lap_{D}}$ at $s=T$,  the correction $z_g$ has the structure
of an improper Bochner integral   converging in $L_2(\Omega)$, namely 
\begin{equation}
  \label{zg-intro}
  z_g = \dashint_0^T \lap e^{(T-s) \lap_{D}} K_0 g(s) \,ds.  
\end{equation}
Hereby the Poisson operator $K_0\colon H^{1/2}(\partial\Omega)\to Z(\mlap)$ is chosen as the inverse of the operator, which results by restricting the boundary trace $\gamma_0\colon H^1(\Omega)\to H^{1/2}(\partial\Omega)$ to its coimage $Z(\mlap)$ of harmonic functions in $H^1(\Omega)$; there is a direct sum $H^1(\Omega)=H^1_0(\Omega)\dotplus Z(\mlap)$.

It is noteworthy that the full influence of the boundary data $g$ on the final state $u(T)$ is given
in the formula  for $z_g$ above.
In addition $z_g\colon H^{1/2}(\,]0,T[\,\times\partial\Omega)\to L_2(\Omega)$ is bounded.

\begin{thm}  \label{yz-thm}
For given data $f \in L_2(0,T; H^{-1}(\Omega))$, $g \in H^{1/2}(\,]0,T[\,\times\partial\Omega)$,  
$u_T \in L_2(\Omega)$ 
the final value problem \eqref{heat-intro} is solved by a function $u$ in the Banach space  $X_1$, whereby 
\begin{equation} \label{X1-eq}
\begin{split}
   X_1&= L_2(0,T; H^1(\Omega)) \bigcap C([0,T];L_2(\Omega)) \bigcap H^1(0,T; H^{-1}(\Omega)),
   \\
  \|u\|_{X_1}=&\big(\int_0^T (\|u(t)\|_{H^1(\Omega)}^2   +\|u'(t)\|_{H^{-1}(\Omega)}^2)\,dt
                 +\sup_{0\le t\le T}\|u(t)\|_{L_2(\Omega)}^2\big)^{1/2},
  \end{split}
\end{equation}
if and only if the data in terms of \eqref{yf-intro} and \eqref{zg-intro} satisfy the compatibility
condition 
\begin{equation}  \label{yz-cnd}
  u_T - y_f + z_g \in D(e^{-T \lap_{D}}).
\end{equation}
In the affirmative case, $u$ is uniquely determined in $X_1$ and has the representation
\begin{align}\label{yz-id}
  u(t) = 
   e^{t \lap_{D}} e^{-T \lap_{D}} (u_T - y_f + z_g) 
   + \int_0^t e^{(t-s)\lap} f(s) \,ds - \dashint_0^t \lap e^{(t-s) \lap_{D}} K_0 g(s) \,ds,
\end{align}
where the three terms all belong to $X_1$ as functions of $t$. 
\end{thm}

Clearly the space of admissible data $Y_1$ is here a specific subspace of 
\begin{equation} \label{Y1-eq}
  L_2(0,T;H^{-1}(\Omega))\oplus H^{1/2}(\,]0,T[\,\times\partial\Omega) \oplus L_2(\Omega),
\end{equation}
for by setting $\Phi_1(f,g,u_T)= u_T-y_f+z_g$ we have 
\begin{align} \label{Y1-id}
  Y_1= \left\{(f,g,u_T) \mid u_T - y_f + z_g \in D(e^{-T \lap_{D}}) \right\}
     =D(e^{-T \lap_{D}} \Phi_1).
\end{align}
Here $e^{-T \lap_{D}} \Phi_1$ is an unbounded operator from the space in \eqref{Y1-eq} to $H$.
Therefore $Y_1$ is a hilbertable Banach space when endowed with the corresponding graph norm 
\begin{multline}  \label{Y1-nrm}
  \| (f,g,u_T)\|_{Y_1}^2 =  \|u_T\|^2_{L_2(\Omega)} + \|g\|^2_{H^{1/2}(\,]0,T[\,\times\partial\Omega)} +
  \|f\|^2_{L_2(0,T;H^{-1}(\Omega))}
\\
   + \int_\Omega\Big|e^{-T\lap_{D}}\big(u_T - \int_0^T\!e^{-(T-s)\lap}f(s)\,ds + 
   \dashint_0^{T}\! \lap e^{(T-s)\lap_{D}} K_0 g(s) \,ds\big)\Big|^2\,dx.
\end{multline} 
Using this the solution operator $(f,g,u_T)\mapsto u$ is bounded $Y_1\to X_1$, that is,
\begin{equation}
    \| u\|_{X_1} \le c \| (f,g,u_T) \|_{Y_1}.
 \end{equation}
This can be shown by exploiting the bijection $u(0)\longleftrightarrow u(T)$ to invoke the classical estimates of the initial value problem, which in the present low regularity setting has no compatibility conditions and therefore allows a reduction to the case $g=0$.
So in combination with Theorem~\ref{yz-thm} we have

\begin{thm}
  The final value Dirichlet heat problem \eqref{heat-intro} is well-posed in the spaces $X_1$ and $Y_1$; cf.\ \eqref{X1-eq} and \eqref{Y1-eq}--\eqref{Y1-nrm}.
\end{thm}

The full proofs of the results in this note can be found in our exposition \cite{ChJo17}.



\providecommand{\bysame}{\leavevmode\hbox to3em{\hrulefill}\thinspace}
\providecommand{\MR}{\relax\ifhmode\unskip\space\fi MR }
\providecommand{\MRhref}[2]{%
  \href{http://www.ams.org/mathscinet-getitem?mr=#1}{#2}
}
\providecommand{\href}[2]{#2}

\end{document}